\documentclass[11pt,fleqn,twoside]{article}
\usepackage{amsfonts,amssymb,latexsym}
\makeatletter
\newcommand{\prava}[1]{\small\it
\begin{flushleft}
Copyright \copyright \ 1999 by  #1
\end{flushleft}}

\newcommand{\name}[1]{\begin{flushleft}
                       \LARGE \bf #1
                       \end{flushleft}\vspace{-3mm}}

\newcommand{\Author}[1]{\begin{flushleft}
                       \it #1 \end{flushleft}}

\newcommand{\Adress}[1]{\begin{flushleft}
                       \it #1 \end{flushleft}}

\newcommand{\Date}[1]{\begin{flushleft}
                      \small  \it #1 \end{flushleft}}

\newcommand{\ehkol}{Author \ name}
\newcommand{\ohkol}{Article \ name}
\renewcommand{\@evenhead}{
\hspace*{-3pt}\raisebox{-15pt}[\headheight][0pt]{\vbox{\hbox to \textwidth 
{\thepage \hfil \ehkol}\vskip4pt \hrule}}}
\renewcommand{\@oddhead}{
\hspace*{-3pt}\raisebox{-15pt}[\headheight][0pt]{\vbox{\hbox to \textwidth 
{\ohkol \hfil \thepage}\vskip4pt\hrule}}}
\renewcommand{\@evenfoot}{}
\renewcommand{\@oddfoot}{}

\newcommand{\be}{\begin{equation}}
\newcommand{\ee}{\end{equation}}
\newcommand{\ba}{\hspace*{-5pt}\begin{array}}
\newcommand{\ea}{\end{array}}
\newcommand{\p}{\partial}
\newcommand{\ds}{\displaystyle}
\makeatother

\newtheorem{Def}{Def\/inition}
\newtheorem{theo}{Theorem}

\begin{document}

\thispagestyle{empty}
\setcounter{page}{127}

\renewcommand{\ehkol}{A.V. Shapovalov and A.Yu. Trifonov}
\renewcommand{\ohkol}{Semiclassical Solutions of the Nonlinear Schr\"odinger
Equation}

\begin{flushleft}
\footnotesize \sf Journal of Nonlinear Mathematical Physics \qquad
1999, V.6, N~2, \pageref{shapovalov-fp}--\pageref{shapovalov-lp}.
\hfill {\sc Letter}
\end{flushleft}

\vspace{-5mm}

{\renewcommand{\footnoterule}{} {\renewcommand{\thefootnote}{}
\footnote{\prava{A.V. Shapovalov and A.Yu. Trifonov}}}

\name{Semiclassical Solutions of the Nonlinear Schr\"odinger
Equation}\label{shapovalov-fp}

\Author{A.V. SHAPOVALOV~$^\dag$ and A.Yu. TRIFONOV~$^\ddag$}

\Adress{$\dag$~Department of Theor. Physics, Tomsk State
University, Tomsk 634050, Russia\\ ~~E-mail:
shpv@phys.tsu.ru\\[1mm]
 $\ddag$~Department of Math. Physics,
 Tomsk Polytechnical University, Tomsk 634034, Russia\\
 ~~E-mail: trifonov@phtd.tpu.edu.ru }

\Date{Received November 7, 1998;  Accepted January 20, 1999}

\begin{abstract}
\noindent
 A concept of semiclassically concentrated solutions is
formulated for the multidimensional  nonlinear Schr\"odinger
equation (NLSE) with an external f\/ield. These solutions are
considered as multidimensional solitary waves. The center of mass
of such a solution is shown to move along with the
bi\-charac\-teris\-tics of the basic symbol of the corresponding
linear Schr\"o\-din\-ger equation. The leading term of the
asymptotic WKB-solution is constructed for the multidimensional
NLSE. Special cases are considered for the standard
one-dimensional NLSE and for NLSE in cylindrical coordinates.
\end{abstract}

\section{Introduction}\label{shapovalov:Intr}

Soliton phenomena is an attractive f\/ield of present day research
in nonlinear physics and mathematics. Essential ingredients in the
soliton theory are the nonlinear Schr\"odinger equation (NLSE) and
its variants appearing in a wide spectrum of problems. Examples
are coupled non\-li\-near
optics~\cite{shapo:Hasegawa,shapo:Mollenauer,shapo:Gordon,shapo:Nassar},
superconductivity~\cite{shapo:Tinkham,shapo:Floria}, and excitation in lattice
systems~\cite{shapo:Christiansen}.

More exactly, solitons are identif\/ied with a certain class
of ref\/lectionless solutions of the equations integrable via the
Inverse Scattering Transform  method (see, for
example,~\cite{shapo:Zakh1}). Such equations, including NLSE, are named
soliton equations. At every instant a soliton is localized in a
restricted spatial region with its centroid moving like a
particle. The particle-like properties of solitons are also
manifested in their elastic collisions.

Soliton equations make up a narrow class of nonlinear equations,
whereas a wider set of nonlinear equations, being nonintegrable
in the framework of the IST, possess soliton- like solutions.
They are localized in some sense, propagate with small energy losses,
and collide with a varied extent of inelasticity.
These solutions are termed soli\-ta\-ry waves (SWs),
quasisolitons, soliton-like solutions, etc. to dif\/ferentiate them
from the solitons in the above exact meaning.
The stability of the localized form of solitons and SWs and their
elastic collisions have led to interesting  physical applications.

It is of interest to study the inf\/luence of  external f\/ields on
the soliton propagation. To do that it is necessary to modify the
original solitary equation by introducing variable coef\/f\/icients
representing an external f\/ield potential that breaks the
IST-integra\-bi\-lity. This problem was considered by
variatio\-nal methods~\cite{shapo:Anderson}, the theory of
soliton perturbations (see, for example,
reviews~\cite{shapo:Kivshar,shapo:Maymistov}), or by using an appropriate
ansatz~\cite{shapo:Nassar,shapo:Moura}.

An investigation of soliton-like states is a separate problem for
multi\-di\-men\-sional mo\-dels. The IST constructions are
unsuitable for them except for the D4 self-dual Yang-Mills
equations and their reductions in $1\le D \le 3$ (see, for
example,~\cite{shapo:Ivanova} for details). The problem is more
complicated in view of the singular behavior of NLSE in
2D~\cite{shapo:Zakh2}. Nevertheless, soliton-like solutions for NLSE on
a plane can exist, as shown in~\cite{shapo:Jackiw} in terms of a
suitable ansatz.

Since the methods for constructing exact
solutions for multidimensional models are restricted if compared
to the cases of 2D models, approximate methods
should be used. An ef\/fective approach to the problem can be
developed based on the WKB-method~\cite{shapo:Maslov}.

}

Consider the generalized nonlinear Schr\"odinger equation for a
``matter'' f\/ield $\Psi (\vec x,t)$:
\be \label{shapovalov:1}
\ba{l}
\ds
-i\hbar \frac {\p}{\p t} \Psi+ \hat {\mathcal H}_{\rm nl}(\Psi )
  = \Biggl \{ -i\hbar \frac {\p}{\p t} + \frac
{1}{2m}\left(-i\hbar \nabla -\vec {\mathcal A}(\vec x,t)\right)^2
\vspace{3mm}\\
\ds \hspace*{50mm} + V (\vec x,t)
-2r|\Psi (\vec x,t)|^2 \Biggr\}\Psi (\vec x,t)=0.
\ea
\ee
Here,  $\vec x \in {\mathbb R }^n$, $t\in {\mathbb R} ^1$; $V (\vec x,t) $,
$\vec{\mathcal A}(\vec x,t)$ are given functions;
$m$ is a real constant, $r$ is a real  parameter of nonlinearity;
$\hbar$ is Planck's constant playing the role of
an asymptotic parameter; $|\Psi|^2=\Psi^*\Psi$,
$\Psi^*$ is the complex conjugate of $\Psi$.

In such a form  equation~(\ref{shapovalov:1}) was
studied in Refs.~\cite{shapo:Tinkham,shapo:Jackiw}.
Since the f\/ield of application of~(\ref{shapovalov:1}) is wider,
the physical meaning of the quantities entering in (\ref{shapovalov:1})
may be quite dif\/ferent from the quantum-mechanical one.
In particular, in the one-dimensional problem of the propagation
of an optical pulse~\cite{shapo:Mollenauer,shapo:Gordon} we have to assume
for  equation~(\ref{shapovalov:1}) $\vec x=x \in {\mathbb R} ^1$ with~$x$
being a normed temporal variable  and $t$ being a normed
space coordinate along which the pulse propagates.
The function $\Psi (x,t)$ is an envelope of the pulse f\/ield.
Equation (\ref{shapovalov:1}) takes the form
\be \label{shapovalov:2}
\left \{-i\hbar \frac {\p}{\p t} +
\frac {1}{2m}\left(i\hbar  \frac {\p }{\p x} + {\mathcal A}(x,t)\right)^2 + V(x,t) -
 2r|\Psi (x,t)|^2 \right\}\Psi (x,t)=0.
\ee
The functions $V(x,t)$ and ${\mathcal A}(x,t)$ simulate the heterogeneity of
the me\-dium. The symmetry of equation~(\ref{shapovalov:2}) was considered in
Ref.~\cite{shapo:Broadbridge}, and asymp\-to\-ti\-cal solutions where
studied  for~(\ref{shapovalov:2}) with ${\mathcal A} =0$, $r=r(x,t)$
in Ref.~\cite{shapo:Subochev}.

In the present work we have formulated the concept of semiclassi\-cally
localized solutions  for~ (\ref{shapovalov:1}), following the
ideas of Ref.~\cite{shapo:Trifonov}. These solutions are the
multidimensional analogues of the soliton, like the solutions
for~(\ref{shapovalov:1}).
The particle-like properties of a solitary wave are
described in terms of the wave centroid. The latter
is shown to move along with the bi\-charac\-teris\-tics of
the basic symbol of the correspondent {\it linear}
Schr\"o\-din\-ger equation.
We construct asympto\-ti\-cal WKB-solutions for equation~(\ref{shapovalov:1})
and consider some examples.

\renewcommand{\thefootnote}{\arabic{footnote}}
\setcounter{footnote}{0}

\section{Semiclassical concentraited solutions} \label{shapovalov:SCS}

Soliton solutions are known to show particle properties.
In classical mechanics, a particle is completely described by its phase
orbit. Therefore, it is natural to introduce a similar concept
for the soliton-like solutions of the nonlinear  Schr\"odinger equation~(\ref{shapovalov:1}).
The way by which to introduce the phase orbit seems
to be obvious enough. It is based on the fact that in quantum mechanics,
the f\/irst moments of a state $\Psi (\vec x,t)$ play the role
of the phase orbit for the quantum system.
Let now $\Psi (\vec x,t)$ be a solution of NLSE~(\ref{shapovalov:1}).
The generalized position operators are $\hat{\vec x} (=\vec x)$
and their conjugate momentum variables are $\hat{\vec p}(=-i\hbar \nabla)$,
\[
[{\hat x}_k, {\hat p}_s] = i\hbar \delta_{k,s},
\qquad k,s=\overline{1,n}.
\]
The mean value of an operator $\hat A$ by the solution $\Psi$
is def\/ined as
\begin{equation} \label{shapovalov:3}
\langle A \rangle =\frac{\langle\Psi|{\hat A}|\Psi\rangle}{\|\Psi\|^2}.
\end{equation}
Here,
${\|\Psi\|^2}=\langle\Psi|\Psi\rangle $;
\[
\langle\Psi|{\hat A(t)}|\Psi\rangle =
\int \Psi^*(\vec x,t){\hat A(t)}\Psi (\vec x,t) d\vec x,
\]
is a function of time for every operator $\hat A(t)$ and it parametrically
depends on $\hbar$,
\begin{equation} \label{shapovalov:4}
\langle \vec x \,\rangle = \vec x(t,\hbar ), \qquad
\langle \vec  p \,\rangle = \vec p(t,\hbar ).
\end{equation}
If there exist the limits:
\begin{equation} \label{shapovalov:5}
\lim_{\hbar \to 0}\vec x(t,\hbar ) =\vec x(t), \qquad
\lim_{\hbar  \to 0}\vec p(t,\hbar ) =\vec p(t),
\end{equation}
then $\vec x(t)$ and $\vec p(t)$  are natural to be named
the phase orbit of the classical system corresponding to the
given solution $\Psi$. It is obvious that both the mean values~(\ref{shapovalov:4})
and the limit values~ (\ref{shapovalov:5}) depend on
the solution $\Psi$ in the general case. Hence, the choice of the
solution  $\Psi$ meets the requirement for the expressions~(\ref{shapovalov:5})
to be a solution of the classical equations of motion.
By analogy with quantum mechanics (see  Ref.~\cite{shapo:Trifonov}),
we can def\/ine  the soliton-like
solutions asymptotic in $\hbar  \to 0$  as follows.

\begin{Def}\label{shapovalov:D1}
Let $z(t)=\{(\vec x(t),\vec p(t)),\,0\leqslant t\leqslant T\}$
be an arbitrary phase orbit in ${\mathbb R}^{2n}$.
We name the solution  $\Psi(\vec x,t,\hbar )$ of equation~(\ref{shapovalov:1})
as a semiclassically concentrated solution ({\sl SCS}) of the class
${\mathbb C}{\mathbb S}_S(z(t),N)$ $(\Psi\in{\mathbb C}{\mathbb S}_S(z(t),N))$
if:

{\rm (i)} there exist the generalized limits
\footnote{By generalized limit we mean the passage to the limit
standardly defined in the distribution theory
(see, for example,  Ref. \cite{shapo:Rudin}).}
\[
 \lim_{\hbar \to0}\frac{|\Psi(\vec x,t,\hbar )|^2}{\|\Psi\|^2}=\delta (\vec x-\vec x(t)),
\qquad
\lim_{\hbar \to0}\frac{|\tilde\Psi(\vec p,t,\hbar )|^2}{\|\Psi\|^2}=\delta (\vec p-\vec p(t));
\]

{\rm (ii)} there exist the centered moments
\[
\Delta ^{(k)}_{\alpha,\beta}(t,\hbar ) =
\langle \hat \Delta^{(k)}_{\alpha,\beta} \rangle,
\qquad 0\leqslant k\leqslant N.
\]
\end{Def}

Here, $\alpha$ and $\beta$ are multiindices, $|\alpha|+|\beta|=k$,
$0\leqslant k\leqslant N$, and
$\hat\Delta^{(k)}_{\alpha,\beta}$ is an operator with
the symmetrized (Weyl) symbol
$\Delta^{(k)}_{\alpha,\beta}(\vec p,\vec x)=$
$(\vec p-\vec p(t))^\alpha (\vec x-\vec x(t))^\beta$.
Recall that
the multiindex $\alpha$ is a vector of the form
$\alpha=(\alpha_1,\dots,\alpha_n)$, where
$\alpha_j\geqslant0$ are integer numbers.
In addition,
$|\alpha|=\sum\limits_{j=1}^n\alpha_j$
and for a vector $\vec\zeta=(\zeta_1,\dots,\zeta_n)
\in{\mathbb R}^n$ we suppose
$\vec\zeta^\alpha =\prod\limits_{j=1}^n
\zeta_j^{\alpha_j}$.

Note that the vectors $\vec x(t)$ and $\vec p(t)$ in
Def\/inition~\ref{shapovalov:D1}  are by no means connected with each other.
The vector $z(t)=(\vec p(t),\vec x(t))$ is named
the classical phase orbit of the system.

Def\/inition~\ref{shapovalov:D1} specif\/ies the concept of solitary waves
for asymptotic solutions of NLSE~(\ref{shapovalov:1}).

\begin{theo}
If a solution  $\Psi$ of  (\ref{shapovalov:1}) is semiclassically
concentrated $(\Psi$ $\in $ ${\mathbb C}{\mathbb S}_S(z(t),N))$,
then $z(t)=(\vec p(t)$, $\vec x(t))$ is a solution of
a classical Hamilton system with the Hamiltonian
\[
{\cal H}_{\rm cl}(\vec p, \vec x, t) =
\frac {1}{2m}\left(\vec p -\vec{\mathcal A}(\vec x,t)\right)^2 + V (\vec x,t).
\]
\end{theo}

\noindent
{\bf Proof}. Let $\hat A$ be an operator, then
the Ehrenfest theorem~\cite{shapo:Ehrenfest}
is true for the mean value of $\hat A$:
\[
\frac {d}{d t} \langle A \rangle =\langle \frac {\p \hat A}{\p t} \rangle +
\frac {i}{\hbar } \langle [\hat{\mathcal H}_{nl} ,\hat A]\rangle .
\]
In particular, for the operators $\hat {\vec p}$ and $\hat {\vec  x}$ we have:
\begin{equation}\label{shapovalov:6}
\frac {d}{d t} \langle {\vec p}\,\rangle =
\frac {i}{\hbar } \langle [\hat{\mathcal H}_{\rm nl} ,\hat {\vec
p}\,] \rangle,\qquad
\frac {d}{d t} \langle {\vec x} \,\rangle =
\frac {i}{\hbar } \langle [\hat{\mathcal H}_{\rm nl} ,\hat {\vec
x}\,] \rangle .
\end{equation}
Using the obvious relations
\[
\langle [ |\Psi (\vec x,t)|^2, \hat{\vec p}\,]\rangle =
\langle [ |\Psi (\vec x,t)|^2, \hat{\vec x}\,]\rangle =0,
\]
we have from (\ref{shapovalov:6}):
\begin{equation}\label{shapovalov:7}
\frac {d}{d t} \langle {\vec p}\,\rangle =
\frac {i}{\hbar} \langle [\hat{\mathcal H}_{\rm l} ,\hat {\vec
p}\,]\rangle ,\qquad
\frac {d}{d t} \langle {\vec x} \,\rangle =
\frac {i}{\hbar} \langle [\hat{\mathcal H}_{\rm l} ,\hat {\vec x}\,]\rangle,
\end{equation}
where
\[
\hat{\mathcal H}_{\rm l}=
\frac {1}{2m}\left(i\hbar \nabla  +\vec {\mathcal A}(\vec x,t)\right)^2 +V(\vec x,t).
\]
With Def\/inition~\ref{shapovalov:D1} of the SCS further proof coincides with a similar one
for the linear case~\cite{shapo:Trifonov}.

\medskip

\noindent
{\bf  Remark.}  Emphasize that, as follows from the theorem,
the centroid of the SCS $\Psi$  moves along the bicharacteristics of
the linear Schr\"odinger equation.

\section{Asymptotic solutions } \label{shapovalov:Sol}

In Section~\ref{shapovalov:SCS} we have discussed the def\/inition and
the basic features of SCS.
Here, we study a theoretical possibility of construction of
the WKB-asymptotic semiclassically concentrated
solutions of equation~(\ref{shapovalov:1}) on a limited time domain $0<t<T$
with  $\hbar$-independent~$T$.

Taking into account the form of the one-soliton solution of
NLSE (see, for example,~\cite{shapo:Zakh1}), let us try solution
of equation~(\ref{shapovalov:1}) in the form
\begin{equation}\label{shapovalov:8}
\Psi =\rho (\theta, \vec x,t,\hbar )
\exp \left[\frac {i}{\hbar }S(\vec x,t,\hbar )\right].
\end{equation}
Here, $\theta =\hbar ^{-1}\sigma(\vec x,t,\hbar )$ is a ``fast'' variable;
$\sigma(\vec x,t,\hbar )$, $\rho (\theta, \vec x,t,\hbar )$, and $S(\vec x,t,\hbar)$
are real functions regular in $\hbar$, that is:
\[
S(\vec x,t,\hbar)= S(\vec x,t)+\hbar S_1(\vec x,t)+\cdots .
\]
The solution (\ref{shapovalov:8}) is assumed to be localized according to Def\/inition~(\ref{shapovalov:D1}).

The derivative operators $\p/\p t$ and $\nabla $ act on the function
(\ref{shapovalov:8}) as follows:
\[
-i\hbar \frac \p{\p t}= -i\hbar \frac\p{\p t}\Big|_{\theta ={\rm const}}-
i\sigma_{,t} \frac\p{\p \theta},\qquad
-i\hbar\nabla= -i\hbar\nabla|_{\theta ={\rm const}}-
i(\nabla \sigma ) \frac\p{\p \theta},
\]
where $\sigma _{,t}=\p \sigma/\p t$.

Henceforth we put
${\p }/{\p t}|_{\theta ={\rm const}} \equiv \p _t$,
$\nabla  |_{\theta ={\rm const}}\equiv \nabla $,
$\p /\p \theta \equiv \p _\theta$.

Substituting (\ref{shapovalov:8}) into (\ref{shapovalov:1}) we f\/ind:
\[
\ba{l}
\ds \exp \left(\frac {i}{\hbar}S\right)\Biggl\{-i\hbar {\p}_t
+ S_{,t} - i\sigma _{,t} \p _\theta + V
- \frac {\hbar ^2}{2m}\nabla ^2 -
\frac {\hbar }{2m}(\nabla ^2\sigma )\p _\theta -\frac {\hbar }{m}(\nabla \sigma \cdot \nabla )\p _\theta
\vspace{3mm}\\
\ds  \quad - i\frac {\hbar }{2m}(\nabla ^2 S) -
\frac {1}{2m}(\nabla \sigma )^2 {\p ^2}_{\theta \theta} -
i\frac {\hbar}{m}(\nabla S \cdot \nabla )-
i\frac {1 }{m}(\nabla S\cdot \nabla \sigma )\p _\theta+
\frac {1}{2m}(\nabla S)^2
\vspace{3mm}\\
\ds \quad + i\frac {\hbar}{m}({\vec A}\cdot \nabla ) +
i\frac {1}{m}({\vec A}\cdot \nabla \sigma )\p _\theta -
 \frac {1}{m}({\vec A}\cdot \nabla S) +
i\frac {\hbar }{2m}(\nabla {\vec A}) + \frac {1}{2m}{\vec A}^2 -
2r \rho ^2  \Biggr\}\rho =0.
\ea
\]

Let us gather $\hbar $-free terms in this equation and put their sum to zero.
In the obtained equation we separate real and imaginary parts
and then separate the ``fast'' variable $\theta$ from others.
As a result we come to the following system of equations
which determines the leading term of
the asymptotic solution:
\begin{equation}\label{shapovalov:9}
\sigma _{,t} +\frac 1m\langle\left(\nabla S-\vec{\mathcal A}\right),\nabla \sigma \rangle=0,
\end{equation}
\begin{equation}\label{shapovalov:10}
S _{,t} + V+ \frac 1{2m}\left(\nabla S-\vec{\mathcal A}\right)^2=r \tilde b(t,\vec x),
\end{equation}
\begin{equation}\label{shapovalov:11}
\frac {1}{2m}(\nabla \sigma)^2 \rho _{,\theta \theta}
+2 \rho^3=r\tilde b\rho.
\end{equation}

Here, $\langle \vec a, \vec b\rangle $ denotes the Euclidean scalar
product of the vectors:  $\sum\limits^n _{j=1} a_j b_j$;
the function $\tilde b(t,\vec x)$ appears as a ``separation parameter''
in separating the ``fast'' variable $\theta$, and
$\tilde b(t,\vec x)$ is to be determined in what follows.

Let us look for $\rho (\theta ,\vec x,t, \hbar )$ in the class of functions
satisfying the conditions:
\begin{equation}\label{shapovalov:12}
\lim_{\theta \to \infty }\rho (\theta , \vec x, t, \hbar)=
\lim_{\theta \to \infty }\rho _{,\theta } (\theta , \vec x, t, \hbar)=0.
\end{equation}
Integrating equation~(\ref{shapovalov:11}) in view of (\ref{shapovalov:12}) we obtain:
\begin{equation}\label{shapovalov:13}
\rho _{,\theta}=
\sqrt {\frac {2m r}{(\nabla \sigma )^2}}
\sqrt{\tilde b -\rho^2} \rho.
\end{equation}
Let us put $r=\varkappa^2>0$ that corresponds to the existence
of soliton solutions in the case when equation~(\ref{shapovalov:1})
is reduced to the standard NLSE (one-dimensional with  $V (\vec x,t) =$,
$\vec{\mathcal A}(\vec x,t)=0$)~\cite{shapo:Zakh1}.
Then $\tilde b=b^2>0$ and further integration of~(\ref{shapovalov:13}) results in
\begin{equation}\label{shapovalov:14}
\rho =\frac {b}{\cosh \left\{b\varkappa\sqrt {2m (\nabla \sigma)^{-2}}
\left(\frac {1}{\hbar }\sigma (\vec x,t)+\sigma _1(\vec x,t)\right)\right\}},
\end{equation}
where the function $\sigma _1(\vec x,t)$ appears as a
``separation constant'' which is to be determined later.

By def\/inition, the ``fast'' variable $\theta$  must have the structu\-re $\frac {1}{\hbar} \sigma (\vec x,t)$.
For equation~(\ref{shapovalov:14}) to correspond this condition it is necessary to set
\[
b\varkappa \sqrt {\frac {2m}{(\nabla \sigma )^2}}=
{\rm const}=\alpha.
\]
Without loss of generality we can put  $\alpha =1$ and then we have
\begin{equation}\label{shapovalov:15}
\rho = \sqrt {\frac {(\nabla \sigma )^2}{2m \varkappa ^2}} \;
\frac {1}{\cosh \left(\frac {1}{\hbar }\sigma (\vec x,t)+\sigma _1(\vec x,t)\right)}.
\end{equation}
Equation (\ref{shapovalov:10}) takes the form
\[
S _{,t} + V+ \frac 1{2m}\left(\nabla S-\vec{\mathcal A}\right)^2=\frac {1}{2m}(\nabla \sigma)^2,
\]
that, together with (\ref{shapovalov:9}), is equivalent to the single
comp\-lex Ha\-mil\-ton-Jaco\-bi equat\-ion
\begin{equation}\label{shapovalov:16}
(S+i\sigma )_{,t} + V + \frac {1}{2m}\left[\nabla (S+i\sigma )
-\vec{\mathcal A}\right]^2=0.
\end{equation}
Thus, for the leading term of the asymptotic expansion~(\ref{shapovalov:8}),
\[
\Psi (\vec x,t,\hbar )=\Psi ^0(\vec x,t,\hbar ) + O(\hbar ),
\]
we have:
\begin{equation}\label{shapovalov:17}
\Psi ^0 =\rho (\theta, \vec x,t,\hbar)
\exp \left[\frac {i}{\hbar}S(\vec x,t)+S_1(\vec x,t)\right],
\end{equation}
where $\rho$ has the form (\ref{shapovalov:15}) and the functions
$\sigma _1(\vec x,t), S_1(\vec x,t)$ are determined from
successive approximations.
The function~(\ref{shapovalov:17}) can be  represented in the form:
\be \label{shapovalov:18}
\Psi (\vec x, t, \hbar)=
2 \sqrt {\frac {(\nabla \sigma )^2}{2m \varkappa ^2}}\;
\frac {\Psi _0(\vec x, t, \hbar)}{1+
|\Psi _0(\vec x, t, \hbar)|^2},
\ee
where
\[
\Psi _0(\vec x, t, \hbar)=\exp \left \{\frac {i}{\hbar} \left[S(\vec x,t)+
i\sigma (\vec x,t)+ \hbar (S_1(\vec x,t)+i\sigma _1(\vec x,t))\right] \right\}.
\]
It can easily to show that $\Psi _0$ is an asymptotic solution of
the  {\it linear} Schr\"odinger equation:
\be \label{shapovalov:19}
\left \{ -i\hbar \frac {\p}{\p t} + V (\vec x,t)
+\frac {1}{2m}\left(-i\hbar \nabla -\vec{\mathcal A}(\vec x,t)\right)^2 \right\}
\Psi _0 (\vec x,t)=O(\hbar ^\alpha ),
\ee
where $\alpha =1$.
Let us write a function $\Psi$, which satisf\/ies  equation~(\ref{shapovalov:1}) to an accuracy
of~$O(\hbar ^2)$, in the form
\be \label{shapovalov:20}
\Psi = \Psi ^0(1+\hbar \Psi ^1).
\ee
Here, $\Psi ^0$  is determined by expression
(\ref{shapovalov:17}) and the function  $\Psi ^1(\theta ,\vec x, t)$
is to be de\-ter\-mined.
As can be seen from (\ref{shapovalov:15}), it is convenient to take
the variable $\theta$ in the form
\be \label{shapovalov:21}
\theta =
\frac 1\hbar \sigma (\vec x,t)+\sigma _1(\vec x,t).
\ee
Let us denote
${\rm Re}\, \Psi^1(\theta,\vec x,t)=u(\theta , \vec x,t)$,
${\rm Im}\,\Psi ^1(\theta , \vec x,t)=v(\theta , \vec x,t)$.
Substituting (\ref{shapovalov:20}) into~(\ref{shapovalov:1}) and
setting to zero
sum\-mands at the equal powers of $\hbar $.
Extract real and imaginary parts in the obtained equations
and f\/ind a system of equations for function~(\ref{shapovalov:20}).
This system includes the Hamilton-Jacobi equation~(\ref{shapovalov:16}) for the
functions $S(\vec x,t)$, $\sigma (\vec x, t)$
and also the following equations determining the functions
$S _1(\vec x,t)$ and $\sigma _1(\vec x, t)$:
\be
S_{1,t}+\frac 1m\langle \left(\nabla S-\vec{\mathcal A}\right),\nabla S_1 \rangle -
\frac 1m\langle \nabla \sigma,\nabla \sigma _1\rangle +
\frac {1}{2m}\Delta \sigma+
\frac {1}{2m}\langle \nabla \sigma,\nabla \rangle
\log (\nabla \sigma )^2 =0,   \label{shapovalov:22}
\ee
\be
\ba{l}
\ds \sigma _{1,t}+\frac 1m\langle \left(\nabla S-\vec{\mathcal A}\right),
\nabla \sigma_1\rangle+
\frac 1m \langle \nabla \sigma, \nabla S_1\rangle
\vspace{3mm}\\
\ds \qquad -  \frac 12\left [\frac 1m\langle \nabla ,\left(\nabla S -\vec{\mathcal A}\right)\rangle +
\left[\log (\nabla \sigma )^2\right]_{,t}+
\frac 1m\langle \left(\nabla S-\vec{\mathcal A}\right),\nabla \log (\nabla \sigma)^2
\rangle \right]=0.
\ea \hspace{-4.1pt}\label{shapovalov:23}
\ee
The functions $u(\theta , \vec x,t)$, $v(\theta , \vec x,t)$
are given by the expressions
\be
\ba{l}
\ds \rho (\theta ,\vec x,t)u(\theta ,\vec x,t)=
\sqrt { \frac {2m}{\varkappa ^2(\nabla \sigma )^2}}\;
\frac {1}{\cosh \theta } \Biggl\{ C_1(\vec x,t) \tanh \theta+  \frac {1}{2m}\langle \nabla \sigma ,
\nabla \sigma _1\rangle
\vspace{3mm}\\
\ds \qquad + \frac {1}{12m}\left[\Delta \sigma +
\langle \nabla \sigma , \nabla \log (\nabla \sigma )^2\rangle \right]
\left(\sinh \theta \cosh \theta- \epsilon  \cosh ^2 \theta \right)\Biggr\},
\ea \label{shapovalov:24}
\ee
\be
\ba{l}
\ds  \rho (\theta ,\vec x,t)v(\theta ,\vec x,t)=
\sqrt { \frac {2m}{\varkappa ^2(\nabla \sigma )^2}}
\Biggl\{ \frac {C_1(\vec x,t)}{ \cosh \theta }+ \frac {1}{4}\left[\frac 1m \langle \nabla ,\nabla S -
\vec{\mathcal A} \rangle \right.
\vspace{3mm}\\
\ds \left. \qquad + \left(\p _t + \frac {1}{m}\langle \left(\nabla S-\vec{\mathcal A}\right),
\nabla \rangle \right) \log (\nabla \sigma )^2\right]
(\epsilon  \sinh \theta -\cosh \theta ) \Biggr\}.
\ea \label{shapovalov:25}
\ee
Here, $\epsilon ={\rm sign} (\sigma)$ and the function $C_1(\vec x, t)$
is determined by successive approximations.

\medskip

\noindent
{\bf Remark.} From relations (\ref{shapovalov:22}) and (\ref{shapovalov:23}) it follows that (\ref{shapovalov:19})
is accurate to $O(\hbar ^2)$.

\medskip

Relation  (\ref{shapovalov:18}) can be considered  as a transformation
con\-nec\-ting the asymptotic solutions of the nonlinear equation~(\ref{shapovalov:1})
with the solution of the linear equation (\ref{shapovalov:19}) for $\alpha=2$.

So, the leading term $\Psi^0$ of the asymptotic solution
of equation~(\ref{shapovalov:1}) is completely f\/ixed by expressions
(\ref{shapovalov:15}), (\ref{shapovalov:17}), (\ref{shapovalov:20}), and~(\ref{shapovalov:21}),
where the functions $S(\vec x,t)$,
$\sigma (\vec x,t)$, $S_1 (\vec x,t)$, and $\sigma _1(\vec x,t)$
are determined by the system of equations~(\ref{shapovalov:16}),
(\ref{shapovalov:22}), and~(\ref{shapovalov:23}).

\section{Special solutions }\label{shapovalov:Special}

For a more detailed study of the above asymptotic solutions~(\ref{shapovalov:17}),
let us consider some special cases of these solutions.

\subsection{One-dimensional NLSE }\label{shapovalov:NS}

Let us put   $\vec x=x \in {\mathbb R}^1$, $\Delta =\p ^2/\p x^2$ and
$\vec{\mathcal A}=V=0$ in (\ref{shapovalov:1}); then equation~(\ref{shapovalov:1}) takes the form
\be\label{shapovalov:26}
\left[i\hbar \frac{\p}{\p t} + \frac {\hbar ^2}{2m} \frac{\p ^2}{\p x^2}+
2\varkappa ^2 |\Psi (x,t,\hbar)|^2\right]\Psi (x,t,\hbar )=0.
\ee
The leading term of the asymptotic solution (\ref{shapovalov:20})
is as follows:
\be \label{shapovalov:27}
\Psi (x,t)=\rho (x,t,\hbar ) \exp \left[\frac {i}{\hbar }(S(x,t) +
\hbar  S_1(x,t))\right],
\ee
where
\be \label{shapovalov:28}
\rho (x,t,\hbar  )=\sqrt {\frac {\sigma _{,x}^2}{2m \varkappa^2}}
\; \frac {1}{\cosh \left(\frac {1}{\hbar } \sigma (x,t)+
\sigma _1(x,t)\right)}.
\end{equation}
Equations (\ref{shapovalov:16}), (\ref{shapovalov:22}), and (\ref{shapovalov:23}) become:
\be
(S+i \sigma )_{,t}+\frac {1}{2m} (S+i \sigma)_{,x}^2=0,\label{shapovalov:29}
\ee
\be
 S_{1 ,t}+\frac {1}{m} S_{,x}S_{1,x}-\frac 1m\sigma _{,x}\sigma _{1,x}
+\frac {3}{2m}\sigma _{,xx}=0,\label{shapovalov:30}
\ee
\be
 \sigma _{1,t} +\frac {1}{m}S_{,x}\sigma _{1,x}+
\frac {1}{m} \sigma _{,x}S_{1,x}-\frac 12\left[\frac 1m S_{,xx}+
\frac {2}{\sigma _{,x}}\sigma _{,xt}+
\frac {2}{m}\frac {S_{,x}}{\sigma _{,x}}\sigma _{,xx}\right]=0.\label{shapovalov:31}
\ee
To construct solution (\ref{shapovalov:27}) and (\ref{shapovalov:28}),
let us look for a special solution of equation~(\ref{shapovalov:29}) in the form
\be\label{shapovalov:32}
S=\alpha _1t+\alpha _2x +\varphi _0,\qquad
\sigma =\beta _1 t+\beta _2 (x-x_0),
\ee
where $\alpha _1$, $\alpha _2$, $\beta _1$, $\beta _2$, $\varphi _0$,
and $x_0$ are real constants.
Substitution (\ref{shapovalov:32}) into (\ref{shapovalov:29}) gives
\[
\alpha _1= \frac {1}{2m}(\beta _2^2-\alpha _2^2), \quad
\beta _1=-\frac {1}{m}\alpha _2\beta _2.
\]
Substituting (\ref{shapovalov:32}) into (\ref{shapovalov:30}) and (\ref{shapovalov:31}), we easily
obtain the following complex equation for the functions
$S _1$ and $\sigma _1$:
\be \label{shapovalov:33}
(S_1 +i \sigma _1)_{,t}+ \frac {\alpha _2+
i\beta _2}{m}(S_1 +i \sigma _1)_{,x} =0.
\ee
Equation (\ref{shapovalov:33}) is a complex wave equation the solution
of which can be written  as
\be\label{shapovalov:34}
w (x,t)\equiv S_1 (x,t)+i \sigma _1(x,t)=f(x- at).
\ee
Here,  $a\equiv (\alpha _2+ i\beta _2)/m$ and
$f(\zeta )$ are analytical functions of the complex variable $\zeta= x-at$.
Denote $\beta _2=2\eta$ and $\alpha _2=2 \xi$, then solution  (\ref{shapovalov:27})
takes the form
\be \label{shapovalov:35}
\ba{l}
\ds \Psi =- \frac {2\eta }{\varkappa \sqrt 2m}\;
\frac {1}{\cosh \left[\frac {2\eta }{\hbar  }\left(x-x_0-\frac {2\xi}{m}t\right)+
{\rm Im}f(x-at)\right]}
\vspace{3mm}\\
\ds \qquad \times
\exp \left [\frac {i}{\hbar }(2\xi x-\frac {2}{m}
(\xi ^2-\eta ^2)t+\varphi _0+\hbar {\rm Re}\,f(x-at) )\right].
\ea
\ee
If $f (x-at)=0$ in (\ref{shapovalov:34}), then (\ref{shapovalov:35}) takes the form
of the exact one-soliton
solution of the nonlinear Schr\"odinger equation
(\ref{shapovalov:26}) which is reduced to the standard form~\cite{shapo:Zakh1}
when $\hbar =m=1$, $\varkappa ^2=1/2$.

If $f(x-at)\neq 0$, we have an asymptotic solution similar
to the one-soliton solution.

\subsection{NLSE with a separated potential}\label{shapovalov:NSV}

Consider equation~(\ref{shapovalov:26}) with a potential $V(x,t)$:
\be \label{shapovalov:36}
\left[i\hbar  \frac {\p}{\p t} + \frac {\hbar ^2}{2m} \frac{\p^2}{\p x^2} +
2\varkappa ^2 |\Psi (x,t,\hbar )|^2-V(x,t)\right]\Psi (x,t,\hbar )=0.
\ee
Solution (\ref{shapovalov:27}) and (\ref{shapovalov:28}) is determined by
(\ref{shapovalov:30}), (\ref{shapovalov:31}), and
\be \label{shapovalov:37}
(S+i \sigma )_{,t}+\frac {1}{2m} (S+i \sigma )_{,x}^2+V=0.
\ee
For the potential $V$ of the separated form,
\[
V(x,t)=v_0(t)+v_1(x),
\]
we can easily f\/ind for the system (\ref{shapovalov:37}), (\ref{shapovalov:30}), and (\ref{shapovalov:31})
two classes of separated solutions determining~(\ref{shapovalov:27}) and~(\ref{shapovalov:28}).

The f\/irst class is  described by the expressions:
\[
\ba{l}
\ds S(x,t)=c_1 t- \int v_0(t) dt,  \qquad
\sigma (x,t)=\int \sqrt {2m(c_1+v_1(x))} dx,
\vspace{3mm}\\
\ds S_1(x,t)=c_2t+c_3,  \qquad
\sigma _1=(3/2)\log|\sigma_{,x}|+
mc_2 \int \sigma _{,x}^{-1} dx +c_4,
\ea
\]
where $c_1,\dots ,c_4= \mbox{const}$.
The second class is presented by
\[
\ba{l}
\ds S(x,t)=c_3 t- \int v_0(t) dt +c_4,
\qquad \sigma (x,t)=c_1 \left[t-m\int \frac {dx}{p'(x)}\right]+c_2,
\vspace{3mm}\\
\ds S_1(x,t)=a_1t+a_3+f(x),
\qquad \sigma _1=a_2t+a_4+g(x),
\ea
\]
Here,  $c_1,\dots ,c_4; a_1,\dots ,a_4 = \mbox{const}$
and the functions $p(x)$, $f(x)$, $g(x)$ are determined in quadratures
from the following equations:
\[
\ba{l}
\ds \frac 1m (p'(x))^2=-(v_1(x)+c_3)+\sqrt {(v_1(x)+c_3)^2+c_1^2},
\vspace{3mm}\\
\ds g'(x)= \frac {m}{p'(x)}\left [\frac {c_1}{p'(x)} f'(x)- \frac {1}{2m}p''(x)-a_2\right],
\vspace{3mm}\\
\ds  [(p'(x)]^4+c_1^2m^2)f'(x)=c_1m^2a_2p'(x)-ma_1[p'(x)]^3-   c_1mp'(x)p''(x).
\ea
\]
The above solutions show the potentialities of separation
of variables as applied to the solution of the general system
of equations (\ref{shapovalov:16}), (\ref{shapovalov:22}), and (\ref{shapovalov:23}).

\subsection{Cylindrical coordinates}\label{shapovalov:NSC}

The self-focusing ef\/fect of a beam of high power optical radiation
propa\-ga\-ting along the $z$-axis  is described by a
nonlinear Schr\"odinger equation of the form~\cite{shapo:Landau}:
\be \label{shapovalov:38}
\left[i\hbar \frac{\p }{\p t} + \frac {\hbar ^2}{2m} \Delta _{\perp} +
2\varkappa ^2 |\Psi (x,y,t,\hbar )|^2\right]\Psi (x,y,t,\hbar )=0.
\ee
Here, $t$ is the time coordinate in the coordinate system
moving with a group velocity along the direction of
the radiation propagation.
$\Delta _{\perp} =\p ^2/\p x^2+\p ^2/\p y^2$ is the Laplace
operator in a plane orthogonal to the radiation direction.
A function $V(t,x,y)$ simulates some nonstationarity and
heterogeneity of the medium where the pulse is propagating.
In the stationary statement of the problem, the variable $z$
takes the place of $t$.

Let us f\/ind asymptotic solutions of the form (\ref{shapovalov:20})
for equation~(\ref{shapovalov:38}) in the cylindrical coordinates
\be\label{shapovalov:39}
x=r\cos \varphi, \qquad y=r\sin \varphi.
\ee

The leading term  of the asymptotic solution (\ref{shapovalov:20})
takes the form
\be \label{shapovalov:40}
\Psi (r,\varphi ,t,\hbar )= \rho (r,\varphi ,t,\hbar )
\exp \left\{\frac {i}{\hbar}[S(r,\varphi ,t) +
\hbar S_1(r,\varphi ,t)]\right\},
\ee
where
\be \label{shapovalov:41}
\rho (r,\varphi ,t,\hbar )=\sqrt {
\frac {(\nabla _{\perp }\sigma )^2}{2m \varkappa ^2}}
\, \frac {1}{\cosh \left(\frac {1}{\hbar} \sigma (r,\varphi ,t)+
\sigma _1(r,\varphi ,t)\right)},
\ee
$\nabla_{\perp}$ is the gradient operator with respect to
the variables $r$ and $\varphi$ (\ref{shapovalov:39}).

The special solution, illustrating a transverse heterogeneity of
the optical pulse, can be written as follows:
\[
\ba{l}
\ds  \Psi (r,\varphi ,t,\hbar )=\frac {c_1}{\varkappa \sqrt {2m}}
\cosh ^{-1} \left[\left(\frac {c_1}{\hbar }+a_2\right)r+c_1b_1t+
\frac 12 \log r +    \frac {a_1}{\hbar }+a_3\right]
\vspace{3mm}\\
\ds \qquad \times \exp \left \{i \left[\left(\frac{c_1^2}{2m \hbar }+
\frac {a_2c_1}{m}\right)t - mb_1r+\frac {c_2}{\hbar }+c_3\right] \right\},
\ea
\]
where $c_1,c_2,c_3,a_1,a_2,a_3,b_1=\mbox{const}$; $c_1\neq 0$.

\section{Conclusion}\label{shapovalov:Conc}

An outcome of the work are the asymptotic solutions of
the solitary wave type for NLSE~(\ref{shapovalov:1}) obtained for a f\/inite time
interval $0<t<T$.  The question of the validity of long-time asymptotics
requires a special consideration, since, even for the linear case, this
problem has received rather much attention in the literature~\cite{shapo:Voros}.

Solitary wave solutions are considered here similar to quantum
wavepackets. This permits one to investigate the particle-like properties
of the SWs  in terms of the Def\/inition~\ref{shapovalov:D1} based on the  Ehrenfest theorem.

The equations of motion for the centroid of an SW are found to be independent
of the nonlinearity factor $r$ in~(\ref{shapovalov:1}), while the solution itself
depends essentially on $r$.  In such an approach, the centroid of the
solitary wave solution moves like a classical particle in the external
f\/ield described by the potentials $(V, \vec {\mathcal A})$ in~(\ref{shapovalov:1}).
This completely corresponds to the well-known soliton properties
for the one-dimensional case~\cite{shapo:Moura}.
Such a situation takes place only when the nonlinearity
factor is constant. If $r=r(\vec x,t)$,  the f\/ield
$\Psi (\vec x,t,\hbar )$ results in the appearance of additional classical
variables, and the correspondent equations  of motion
become much more complicated.

The problems discussed are beyond the scope of this work and requires a
special investigation.

\subsection*{Acknowledgements}
This work was partially supported by the Russian Foundation
for Basic Research, Grant No 98-02-16195.

\label{shapovalov-lp}
\end{document}